\documentclass[a4paper]{article}
\usepackage[PostScript=dvips]{diagrams}
\usepackage{amssymb, amsthm}
\usepackage{pstricks, pst-node, rotating, url}

\makeatletter
\let\orig@end@rotfloat = \end@rotfloat
\def\end@rotfloat{\orig@end@rotfloat\ignorespacesafterend}
\let\endsidewaysfigure\end@rotfloat
\let\endsidewaystable\end@rotfloat
\makeatother

\psset{nodesep=3pt,linewidth=\fontdimen8\textfont3,%
arrowinset=0, arrowsize=2.8pt 0, arrowlength=1.4, labelsep=3pt}

\newtheorem{theorem}{Theorem}
\newtheorem{propn}[theorem]{Proposition}
\newtheorem{lemma}[theorem]{Lemma}

\theoremstyle{definition}
\newtheorem*{definition}{Definition}

\theoremstyle{remark}
\newtheorem*{remark}{Remark}

\newcommand\C{\mathbb{C}}
\newcommand\D{\mathbb{D}}
\newcommand\tn{\otimes}

\makeatletter
\newenvironment{starredequation}[1][*]{\def\mcr@eqlabel{\textup(#1\textup)}\def\@currentlabel{#1}$$}{\eqno \hbox{\mcr@eqlabel}$$\ignorespacesafterend}
\makeatother

\let\after = \circ

\title{Finite Products are Biproducts in a Compact Closed Category}
\author{Robin Houston\thanks{This work is supported by an EPSRC PhD studentship.}
	\\ \normalsize School of Computer Science, University of Manchester}
\begin{document}
\maketitle

\begin{abstract}\noindent
  If a compact closed category has finite products or finite
  coproducts then it in fact has finite biproducts, and so is
  semi-additive.
\end{abstract}

\section{Introduction}
Compact closed categories with biproducts have recently attracted
renewed attention from theoretical computer scientists, because of
their role in the abstract approach to quantum information initiated by
Abramsky and Coecke \cite{AbramskyCoecke}. Perhaps surprisingly,
it seems to have gone unnoticed that finite products or coproducts in
a compact closed category necessarily carry a biproduct structure.
Here we prove that this is so. In fact we prove a more general result, viz:
{\def\thetheorem{2}\begin{propn}
	Let $\C$ be a monoidal category with finite products and coproducts,
	and suppose that for every object $A \in \C$, the functor $A\tn\mathord-$ preserves
	products and the functor $\mathord-\tn A$ preserves coproducts. Then $\C$ has finite biproducts.
\end{propn}\addtocounter{theorem}{-1}}
A category with finite biproducts is necessarily semi-additive, i.e.\ enriched
over commutative monoids. In other words, each homset has the structure
of a commutative monoid, and composition preserves the commutative monoid
structure. The converse is also true: a semi-additive category with finite
products or coproducts in fact has finite biproducts. Therefore an equivalent
statement of our conclusion would be that $\C$ is semi-additive.

The gist of the argument is as follows. Let $\C$ be as in the statement
of the proposition. The object $0\tn 1$ is initial because $\mathord-\tn1$ preserves
initiality, and terminal because $0\tn\mathord-$ preserves terminality. So it is
a zero object. The binary case is similar, though more intricate. Let $A$, $B$, $C$, $D\in\C$
and consider the object
$(A+B)\tn(C\times D)$. Since the $\tn$ distributes over
both the $+$ and the $\times$ in this expression, it may be multiplied out as either
\begin{equation}\label{mo1}
	(A\tn C \times A\tn D)+(B\tn C \times B\tn D)
\end{equation}
or
\begin{equation}\label{mo2}
	(A\tn C + B\tn C)\times(A\tn D + B\tn D),
\end{equation}
hence (\ref{mo1}) is isomorphic to (\ref{mo2}). Letting $C=D=I$ shows
that
\begin{equation}\label{isom}
	A^2+B^2\cong(A+B)^2,
\end{equation}
and it may be verified (Lemmas~\ref{l1.1}--\ref{l1.2}) that the canonical natural map
\[A^2+B^2\to(A+B)^2,\] denoted $t_{A,B}$ below, is equal to the left-to-right
direction of (\ref{isom}). It follows that $t_{A,B}$ is invertible. From this we
derive, via Lemma~\ref{l1.3}, that the natural map $A+B\to A\times B$ is also invertible,
which implies the desired conclusion.

The remainder of this paper contains the detailed proof. The next section recalls the
basic facts about finite products and coproducts, and some simple properties of compact closed categories: it will not tax the experienced reader, who may prefer to skip directly
to \S\ref{s-main}.


\section{Background}
This short paper uses only elementary ideas of category theory,
which we briefly recall so as to fix our notation.

In a category with finite products, we denote the given terminal object $1$,
and suppose that for every pair $A$, $B$ of objects there is a given
product cone
$(\pi_1: A\times B\to A, \pi_2: A\times B\to B)$.
For any pair of maps $f: X\to A$, $g: Y\to B$, we denote their
pairing as $\langle f,g\rangle: X\to A\times B$, i.e.\ $\langle f,g\rangle$
is the unique map for which
$\pi_1\after\langle f,g\rangle=f$ and $\pi_2\after\langle f,g\rangle=g$.
Given $f: A\to B$ and $g: C\to D$, we write
$f\times g$ for the map
\[
	\langle f\after\pi_1, g\after\pi_2\rangle: A\times C\to B\times D.
\]
Note that this definition makes $\times$ into a functor, in such a
way that $\pi_1$ and $\pi_2$ constitute natural transformations. For
example $\pi_1\after(f\times g) = \pi_1\after\langle f\after\pi_1, g\after\pi_2\rangle
= f\after\pi_1$.

A functor $F$ is said to \emph{preserve products} if the
image under $F$ of a product cone is always a product cone (not necessarily
the chosen one). We
take it to include the nullary case also, i.e.\ the image of
a terminal object must be terminal. If the categories
$\C$ and $\D$ have finite products and $F:\C\to\D$ preserves
products then the morphism
\[
	F(A\times B) \rTo^{\langle F\pi_1, F\pi_2\rangle} FA\times FB
\]
is invertible.

The case of coproducts is dual to the above. In a category that has finite
coproducts, we assume that there is an initial object $0$ and that for
every pair of objects $A$, $B$, there is a given coproduct cocone
$(i_1: A\to A+B, i_2: B\to A+B)$.
Given maps $f: A\to Y$ and $g: B\to Y$, we write their co-pairing as
\[
	[f,g]: A+B\to Y;
\]
if $\C$ and $\D$ have finite coproducts and $F:\C\to\D$ preserves coproducts
then the map
\[
	FA + FB \rTo^{[Fi_1, Fi_2]} F(A+B)
\]
is invertible.

Now suppose we are in a category that has both finite products and
finite coproducts. A morphism
\[
	f: A+B \to C\times D
\]
is determined by the four maps
\[\begin{array}{l@{\;}l@{\qquad}l@{\;}l}
	f_{11} := \pi_1\after f\after i_1:& A\to C,		&		f_{12} := \pi_1\after f\after i_2:& B\to C\\
	f_{21}:= \pi_2\after f\after i_1:& A\to D,		&		f_{22}:= \pi_2\after f\after i_2:&B\to D,
\end{array}\]
since $f = [\langle f_{11}, f_{21} \rangle,\langle f_{12}, f_{22}\rangle]
= \langle[f_{11}, f_{12}], [f_{21}, f_{22}]\rangle$. We refer to this as the
\emph{matrix representation} of $f$, and write it as
\[
	f = \Bigl[\begin{array}[c]{cc}
		f_{11}	&	f_{12}	\\
		f_{21}	&	f_{22}
	\end{array}\Bigr].
\]
A technique that is used several times below is to check that two maps
are equal by calculating and comparing their matrix representations.

There are several equivalent ways of defining what it means for a category to
have finite biproducts. The one most convenient for our purposes is as follows
(see Exercise VIII.2.4 of Mac Lane \cite{MacLane}).
\begin{definition}
	A category $\C$ has \emph{finite biproducts} if it has finite products and
	finite coproducts, such that:
	\begin{itemize}
		\item the unique morphism $0\to 1$ is invertible, thus there is a (unique)
    		zero map $0_{A,B}: A \to 1 \cong 0 \to B$ between any objects $A$ and $B$, and
		\item the morphism
		\[
			\Bigl[\begin{array}{cc}1_A&0_{B,A}\\0_{A,B}&1_B\end{array}\Big]: A+B \to A\times B
		\]
		is invertible for all $A$ and $B$ in $\C$.
	\end{itemize}
\end{definition}

\emph{Compact closed} categories were first defined (almost in passing)
by Kelly \cite{MVFC}, and later studied in depth by Kelly and Laplaza \cite{KL}.
The reader may consult either of those references for the precise definition.
For the purposes of this paper, it suffices to know that a
compact closed category is a monoidal category $(\C,\tn,I)$ that has --
among other things -- the following two properties:
\begin{itemize}
	\item $\C$ is self-dual, i.e.\ $\C$ is equivalent to $\C^{\mathrm op}$,
	\item for every object $A\in\C$, the functors $A\tn\mathord-$ and
		$\mathord-\tn A$ have both a left and a right adjoint.
\end{itemize}
Examples include the category $\mathrm{Rel}$ of sets and relations, with
the tensor as cartesian product, and the category $\mathrm{FinVect}$ of finite-%
dimensional vector spaces, with the usual tensor product of vector spaces.

\section{Main Result}\label{s-main}
Our main result is as follows.
\begin{theorem}\label{theorem}
	Let $\C$ be a compact closed category. If $\C$ has finite products
	(or coproducts) then it has finite biproducts.
\end{theorem}

\noindent We shall deduce the theorem from a somewhat more general proposition:
\begin{propn}\label{prop}
	Let $\C$ be a monoidal category with finite products and coproducts,
	and suppose that for every object $A \in \C$, the functor $A\tn\mathord-$ preserves
	products and the functor $\mathord-\tn A$ preserves coproducts. Then $\C$ has finite biproducts.
\end{propn}

\noindent The nullary case may be dispensed with immediately:

\begin{proof}[Proof that the unique morphism $0\to 1$ is invertible]
	The functor $0\tn\mathord{-}$ preserves products, thus $0\tn 1$ is terminal.
	But also the functor $\mathord-\tn1$ preserves coproducts, so $0\tn 1$ is also initial.
	Therefore $0$ is isomorphic to $1$, and the claim follows.
\end{proof}

From now on, we assume that we have a category that satisfies the
conditions of Proposition~\ref{prop}, and which therefore has a zero object.
We shall omit the subscripts when referring to a zero map, since the
type is always obvious from the context. We have no further occasion
to refer explicitly to an initial object, so the symbol `$0$' below always
denotes a zero map. Also we shall follow the common practice of
abbreviating the identity morphism $1_A$ to $A$.

\begin{remark}
Since $A\tn-$ preserves products, we know that for all
objects $A$,$B$,$C$, the distribution map
\[
  \langle A\tn \pi_1, A\tn \pi_2\rangle: A\tn(B \times C) \to (A\tn B) \times (A\tn C)
\]
is invertible, and since $-\tn C$ preserves coproducts, we know that
for all objects $A$,$B$,$C$, the distribution map
\[
  [i_1\tn C, i_2\tn C]: (A\tn C) + (B\tn C) \to (A + B)\tn C
\]
is invertible.
\end{remark}

\begin{lemma}\label{l1.1}
For all objects $A_1$, $A_2$, $B_1$, $B_2$, the canonical map
\begin{starredequation}[$*$]
   \Bigl[\begin{array}{cc}
    i_1\after\pi_1 & i_2\after\pi_1 \\
    i_1\after\pi_2 & i_2\after\pi_2
   \end{array}\Bigr]
\end{starredequation}
\textup($= [i_1 \times i_1, i_2 \times i_2]
   = \langle \pi_1 + \pi_1, \pi_2 + \pi_2\rangle$\textup)
of type
\[\begin{array}{lcr}
  \begin{array}{l}
  ((A_1\tn B_1) \times (A_1\tn B_2))\\
  \hskip 3em+ ((A_2\tn B_1) \times (A_2\tn B_2))
  \end{array}
  &\!\!\to\!\!&
  \begin{array}{l}
  ((A_1\tn B_1) + (A_2\tn B_1)) \\
  \hskip 3em\times ((A_1\tn B_2) + (A_2\tn B_2))
  \end{array}
\end{array}\]
is invertible.
\end{lemma}
\begin{proof}
We'll show that $(*)$ is equal to the map $y$ defined as the composite
\[\begin{array}{l}
  ((A_1\tn B_1) \times (A_1\tn B_2)) + ((A_2\tn B_1) \times (A_2\tn B_2))\\
  \hskip 3em \to (A_1\tn(B_1\times B_2)) + (A_2\tn(B_1\times B_2))\\
  \hskip 5em \to (A_1 + A_2)\tn(B_1 \times B_2)\\
  \hskip 7em \to ((A_1 + A_2)\tn B_1) \times ((A_1 + A_2)\tn B_2)\\
  \hskip 9em \to ((A_1\tn B_1) + (A_2\tn B_1)) \times ((A_1\tn B_2) + (A_2\tn B_2))
\end{array}\]
of distribution maps and their inverses. Clearly $y$ is invertible,
since it is composed of isomorphisms.

Take $j$, $k \in \{1,2\}$ and consider the diagram in Fig.~\ref{fig-1}.
\begin{sidewaysfigure}
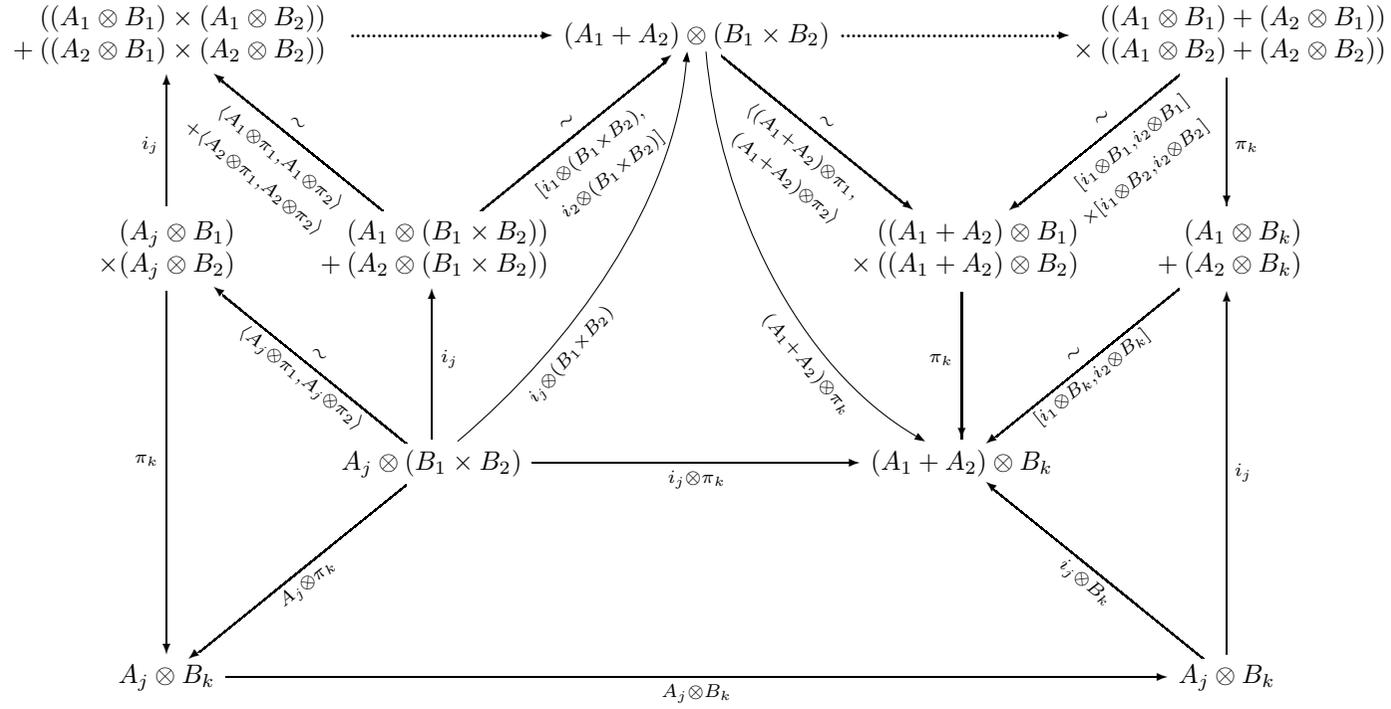

\begin{diagram}[labelstyle=\scriptstyle,w=10em,h=4em,tight]
	\begin{array}{r}((A_1\tn B_1) \times (A_1\tn B_2))\\{}+ ((A_2\tn B_1) \times (A_2\tn B_2))\end{array}
	&\rTo.&
	(A_1 + A_2)\mathbin{\rnode{t}{\tn}}(B_1 \times B_2)
	&\rTo[rightshortfall=-6pt].&
	\begin{array}{r}(
		(A_1\tn B_1) + (A_2\tn B_1))\\
		{}\times((A_1\tn B_2) + (A_2\tn B_2))
	\end{array}
	\\
	\uTo<{i_j}
	&
	\luTo[snake=-5pt](1,2)_{\begin{array}{r}
		\labelstyle\langle A_1\tn\pi_1, A_1\tn\pi_2\rangle\\
		\labelstyle+\langle A_2\tn\pi_1, A_2\tn\pi_2\rangle
	\end{array}}^\sim
	\ruTo(1,2)[snake=0pt]>{\begin{array}{r}
		\labelstyle [i_1\tn(B_1\times B_2),\\
		\labelstyle i_2\tn(B_1\times B_2)]
	\end{array}}^\sim
	& &
	\rdTo(1,2)_{\begin{array}{r}
		\labelstyle \langle(A_1+A_2)\tn\pi_1,\\
		\labelstyle (A_1+A_2)\tn\pi_2\rangle
	\end{array}}^\sim
	\ldTo[snake=8pt](1,2)_{\begin{array}{r}
		\labelstyle [i_1\tn B_1, i_2\tn B_1]\\
		\labelstyle \times[i_1\tn B_2, i_2\tn B_2]
	\end{array}}^\sim
	& \dTo>{\pi_k}
	\\
	\begin{array}{r}
		(A_j\tn B_1) \\
		\times(A_j\tn B_2)
	\end{array}
	&\begin{array}{r}(A_1\tn(B_1 \times B_2))\\{}+(A_2\tn (B_1\times B_2))\end{array}
	&&\begin{array}{r}((A_1+A_2)\tn B_1) \\{}\times((A_1+A_2)\tn B_2)\end{array}
	&\begin{array}{r}(A_1\tn B_k)\\ {}+ (A_2\tn B_k)\end{array}
	\\
	&\luTo(1,2)_{\langle A_j\tn \pi_1, A_j\tn \pi_2\rangle}^\sim
	\rdTo[nohug](0,4)<{\pi_k}
	\uTo>{i_j}
	&&\dTo<{\pi_k}
	\ldTo(1,2)_{[i_1\tn B_k, i_2\tn B_k]}^\sim
	\ldTo[nohug](0,4)<{i_j}
	\\
	&\rnode{l}{A_j\tn(B_1\times B_2)}
	& \rTo_{i_j \tn\pi_k}
	& \rnode{r}{(A_1+A_2)\tn B_k}
	\\
	\ldTo(1,2)>{A_j\tn \pi_k}
	&&&& \luTo(1,2)<{i_j\tn B_k}
	\\
	A_j\tn B_k && \rTo_{A_j\tn B_k} && A_j\tn B_k
	\nccurve[angleA=40, angleB=265, offsetB=3pt]{->}{l}{t}
	\bput{49}(.3){\labelstyle i_j\tn(B_1\times B_2)}
	\nccurve[angleA=275, angleB=150, offsetA=3pt]{->}{t}{r}
	\bput{-49}(.7){\labelstyle (A_1+A_2)\tn\pi_k}
\end{diagram}
\caption{Diagram used in the proof of Lemma~\ref{l1.1}. The arrows marked
	`$\sim$' are invertible, by the remark preceding Lemma~\ref{l1.1}. A dotted
	arrow represents the unique (iso)morphism for which the triangle below it commutes,
	so that the composite along the top edge is equal, by definition, to $y$.}
\label{fig-1}
\end{sidewaysfigure}
All the regions commute for obvious reasons, so the
outside commutes and $\pi_k\after y\after i_j = i_j\after \pi_k$.
Since this is true for all $j$ and $k$, it follows that
$y = (*)$,
as required.
\end{proof}

\begin{definition}
Given objects $A$ and $B$, let $t_{A,B}$ denote the map
\[
    \Bigl[\begin{array}{cc}
    i_1\after\pi_1 & i_2\after\pi_1 \\
    i_1\after\pi_2 & i_2\after\pi_2
   \end{array}\Bigr]: (A\times A) + (B\times B) \to (A+B) \times (A+B)
 \]
 \end{definition}

\begin{lemma}\label{l1.2}
	For all objects $A$, $B$, the map $t_{A,B}$ is invertible.
\end{lemma}
\begin{proof}
Use Lemma~\ref{l1.1} with $A_1=A$, $A_2=B$ and $B_1=B_2=I$, and
apply the right-unit isomorphism.
\end{proof}

\begin{definition}
Given objects $A$ and $B$, let $e_{A,B}$ denote the composite
\[
	(A\times A)+(B\times B) \rTo^{\pi_1 + \pi_2} A+B
	\rTo^{\langle A,0\rangle+\langle 0,B\rangle}
	(A\times A)+(B\times B)
\]
which is clearly an idempotent that splits on $A + B$, and let $e'_{A,B}$
denote the composite
\[
	(A+B)\times(A+B) \rTo^{[A,0]\times[0,B]} A\times B \rTo^{i_1\times i_2} (A+B) \times (A+B)
\]
which is an idempotent that splits on $A\times B$.
\end{definition}

\begin{lemma}\label{l1.3}
$t_{A,B}$ is a map of idempotents from $e_{A,B}$ to $e'_{A,B}$,
i.e.\ the diagram
\begin{diagram}
	(A\times A)+(B\times B) & \rTo^{t_{A,B}} & (A+B)\times(A+B) \\
	\dTo<{e_{A,B}} && \dTo>{e'_{A,B}} \\
	(A\times A)+(B\times B) & \rTo_{t_{A,B}} & (A+B)\times(A+B)
\end{diagram}
commutes.
\end{lemma}
\begin{proof}
We claim that both paths have the matrix representation
\[
	\Bigl[\begin{array}{cc}
		 i_1\after\pi_1 & 0 \\
		 0 & i_2\after\pi_2
	\end{array}\Bigr].
\]
%
Consider the diagram
\begin{diagram}[h=1.5em,w=2em]
	A^2+B^2 & \rTo^{\pi_1 + \pi_2} & A+B
		& \rTo^{\langle A,0\rangle+\langle 0,B\rangle} & A^2+B^2
			& \rTo^{t_{A,B}} & (A+B)^2	\\
		&&&&\uTo>{i_1} \\
		\uTo<{i_1} && \uTo<{i_1} && A^2 && \dTo>{\pi_1} \\
		&&&\ruTo^{\langle A,0\rangle} &\dTo>{\pi_1}\\
		A^2 & \rTo_{\pi_1} &A & \rTo^A & A & \rTo_{i_1} & A+B\\
		&&\dTo>{i_1} & \ruTo_{[A,0]}\\
		\dTo<{i_1} && A+B && \uTo>{\pi_1} && \uTo>{\pi_1}\\
		&&\uTo<{\pi_1}\\
		A^2+B^2 & \rTo_{t_{A,B}} & (A+B)^2
			& \rTo_{[A,0]\times[0,B]} & A\times B & \rTo_{i_1\times i_2} & (A+B)^2
\end{diagram}
where the composite along the top edge is equal to $t_{A,B} \after e_{A,B}$,
and the bottom edge is equal to $e'_{A,B}\after t_{A,B}$. Since all the cells commute, it follows that
\[
	\pi_1\after(t_{A,B}\after e_{A,B})\after i_1  =  i_1\after\pi_1  =  \pi_1\after(e'_{A,B}\after t_{A,B})\after i_1,
\]
and a similar argument shows that
\(
	\pi_2\after (t_{A,B}\after e_{A,B})\after i_2  =  i_2\after \pi_2  =  \pi_2\after (e'_{A,B}\after t_{A,B})\after i_2.
\)
Similar diagrams also show that
\(
	\pi_1\after (t_{A,B}\after e_{A,B})\after i_2  = 0  =  \pi_1\after (e'_{A,B}\after t_{A,B})\after i_2
\)
and
\(
	\pi_2\after (t_{A,B}\after e_{A,B})\after i_1  = 0  =  \pi_2\after (e'_{A,B}\after t_{A,B})\after i_1.
\)
For example, for the former we have
\begin{diagram}[h=1.5em,w=2em]
	A^2+B^2 & \rTo^{\pi_1+\pi_2} & A+B & \rTo^{\langle A,0\rangle+\langle 0,B\rangle}
		&A^2+B^2 & \rTo^{t_{A,B}} & (A+B)^2
	\\
	&\luTo[nohug](0,5)<{i_2}&&&\uTo>{i_2}&\ruTo[nohug](0,5)<{\pi_1}
	\\
	&&\uTo<{i_2}&&B^2
	\\
	&&&\ruTo^{\langle 0,B\rangle}&\dTo>{\pi_1}
	\\
	&&B &\rTo_0&A
	\\
	B^2&\ruTo(2,1)^{\pi_2} &&&&\rdTo(2,1)^{i_2}& A+B
	\\
	&\rdTo[nohug](0,5)<{i_2}\rdTo(2,1)_{\pi_1}&B&\rTo^0&A&\ruTo(2,1)_{i_1}
		\ldTo[nohug,snake=-5pt](0,5)<{\pi_1}
	\\
	&&\dTo>{i_2}&\ruTo_{[A,0]}
	\\
	&&A+B&&\uTo>{\pi_1}
	\\
	&&\uTo>{\pi_1}
	\\
	A^2+B^2&\rTo_{t_{A,B}}&(A+B)^2&\rTo_{[A,0]\times[0,B]}&A\times B
		&\rTo_{i_1\times i_2} & (A+B)^2.
\end{diagram}
\end{proof}

We can now complete the proof of Proposition~\ref{prop}, and hence of Theorem~\ref{theorem}.

\begin{proof}[Proof that $\Bigl[\begin{array}{cc}A&0\\0&B\end{array}{\Bigr]}$ is invertible]
By Lemma~\ref{l1.3}, we know that the map $c_{A,B} :=$
\[
	A+B \rTo^{\langle A,0\rangle + \langle 0,B\rangle}  (A\times A) + (B\times B)
	\rTo^{t_{A,B}} (A+B) \times (A+B) \rTo^{[A,0] \times [0,B]} A\times B
\]
is invertible with inverse
\[
	A\times B \rTo_{i_1\times i_2} (A+B) \times (A+B)
	\rTo_{t_{A,B}^{-1}} (A\times A) + (B\times B)
	\rTo_{\pi_1 + \pi_2} A+B,
\]
so it suffices to check that $c_{A,B} = [\langle A,0\rangle, \langle0,B\rangle]$.
But that's easy to check: for example, the diagram
\begin{diagram}[nohug]
	&&&\rnode{2}{A}
	\\
	\rnode{1}{A} & \rTo_{\langle A,0\rangle} & A^2\ruTo(1,1)_{\pi_1} && A+B\rdTo(1,1)_{i_1}
		& \rTo_{[A,0]} & \rnode{3}{A}
	\\
	\dTo<{i_1} && \dTo>{i_1} &&\uTo<{\pi_1} && \uTo>{\pi_1}
	\\
	A+B & \rTo_{\langle A,0\rangle + \langle0,B\rangle}
		& A^2 + B^2 & \rTo_{t_{A,B}} & (A+B)^2 & \rTo_{[A,0]\times[0,B]}
		& A\times B
	\nccurve[angleA=30,angleB=180]{->}12\Aput{A}
	\nccurve[angleA=0,angleB=150]{->}23\Aput{A}
\end{diagram}
shows that $\pi_1\after c_{A,B}\after i_1$ is the identity on $A$, and the diagram
\begin{diagram}[nohug]
	&&&\rnode{2}{A}
	\\
	\rnode{1}{A} & \rTo_{\langle A,0\rangle} & A^2\ruTo(1,1)_{\pi_2} && A+B\rdTo(1,1)_{i_1}
		& \rTo_{[0,B]} & \rnode{3}{B}
	\\
	\dTo<{i_1} && \dTo>{i_1} &&\uTo<{\pi_2} && \uTo>{\pi_2}
	\\
	A+B & \rTo_{\langle A,0\rangle + \langle0,B\rangle}
		& A^2 + B^2 & \rTo_{t_{A,B}} & (A+B)^2 & \rTo_{[A,0]\times[0,B]}
		& A\times B
	\nccurve[angleA=30,angleB=180]{->}12\Aput{0}
	\nccurve[angleA=0,angleB=150]{->}23\Aput{0}
\end{diagram}
shows that $\pi_2\after c_{A,B}\after i_1 = 0$. The other two cases are similar.
\end{proof}

\begin{proof}[Proof of Theorem~\ref{theorem}]
	A compact closed category is equivalent to its
	opposite, therefore has finite coproducts iff it has finite products. For every
	object $A$, the functors $A\tn-$ and $-\tn A$  have both a left and a right adjoint,
	hence preserve limits and colimits. So Proposition~\ref{prop} applies, in particular,
	to a compact closed category that has finite products (or coproducts).
\end{proof}

\section{Final Remarks}
It is significant that the zero object plays a crucial role in our argument.
A compact closed category may very
well have finite \emph{non-empty} products and coproducts
that are not biproducts. A simple example, due to Masahito
Hasegawa, is the ordered group of integers under
addition. Indeed \emph{any} linearly ordered abelian group constitutes an example,
for the following reason. A partially ordered abelian group may be regarded as a
compact closed category: the underlying partial order is regarded as a category
in the usual way, the group operation provides a symmetric tensor product, and
the adjoint of an object is its group inverse.
If in fact the group is \emph{linearly} ordered then every non-empty finite set
of elements has a minimum (which is their product) and a maximum (coproduct).

This degenerate example may also be used to construct non-degenerate
examples, by taking its product with $\mathrm{Rel}$, say.

One last observation: Proposition~\ref{prop}'s requirement that $\C$ be a monoidal
category is stronger than necessary. We didn't actually need the associativity of tensor,
nor the left unit isomorphism. So instead of the full monoidal structure it suffices merely to
have a functor $\tn: \C \times \C \to \C$ with a right unit.

\section*{Acknowledgements} I am indebted to Peter Selinger for
bringing this question to my attention, and to Robin Cockett for pointing
out how to simplify my original proof.
I have used Paul Taylor's diagrams package.

\bibliography{cs}
\end{document}